\documentclass[11pt,reqno]{amsart}
\usepackage{a4wide}

\usepackage{amssymb, amsmath, amsthm, amsfonts, amscd}
\usepackage[english]{babel}
\usepackage{graphicx}
 \usepackage{amsaddr}
 \usepackage[section]{placeins}

\newtheorem{remark}{Remark}
\newtheorem{thm}{Theorem}

\def\R{\mathbb{R}}
\def\C{\mathbb{C}}

\begin{document}

\title[Neumann controls for the heat equation]{Numerical computation of Neumann  {controls} for the heat equation on a finite interval}
\author{K. Kalimeris, T. Özsarı, and  N. Dikaios}
\thanks{KK and ND were partially supported by the research programme 200/984 of the Research Committee of the Academy of Athens. This work was completed while TÖ was visiting the Academy of Athens. Date of submission: 18/09/2021.}
\thanks{K. Kalimeris is with Mathematics Research Center at Academy of Athens, Athens, 115 27 Greece (e-mail: kkalimeris@academyofathens.gr). }
\thanks{T. Özsarı is with Department of Mathematics, Bilkent University, Ankara 06800 Turkey (e-mail: turker.ozsari@bilkent.edu.tr).}
\thanks{N. Dikaios is with Mathematics Research Center at Academy of Athens, Athens, 115 27 Greece (e-mail: ndikaios@academyofathens.gr).}

\maketitle

\begin{abstract}
	This paper presents a new numerical method which approximates Neumann type null controls for the heat equation and is based on the Fokas method. This is a direct method for solving problems originating from the control theory, which allows the realisation of an efficient numerical algorithm that requires small computational effort for determining the null control with exponentially small error. Furthermore, the unified character of the Fokas method makes the extension of the numerical algorithm to a wide range of other linear PDEs and different type of boundary conditions straightforward.
\end{abstract}


	\section{Introduction}
\label{secintro}
Controllability of evolution equations is an important research topic in engineering and mathematics.  This problem has been studied from many different perspectives.  In physical systems for which access to medium is allowed one can use internal controls acting on the whole domain or only on a subregion within the domain.  In contrast, if access to medium is restricted then one has to work with controls that act externally, namely on the boundary of the domain.  This is generally done in two ways.  The first way is using feedback controllers and construct closed loop systems, with backstepping based boundary controllers being an example, see for instance \cite{KrsticBook} for the design of such control systems.  Feedback controllers are very effective for stabilizing a physical system in long time, preferably with a desired exponential speed.  The second approach for controlling a physical evolution from the boundary of the domain is to design an open loop control system in which the boundary control is simply an explicit function which is fed into the system.  Such controls are generally designed for steering solutions of the physical system to or near a desired target state in finite time.    Exact, null, and approximate controllability systems are some examples of open loop control design.   {It should be pointed out that in practice open loop and closed loop control systems can be used together. A suitable control input can be calculated and fed into a system to steer it near a desired state, while a feedback simultaneously helps to achieve this in a rather fast and uniform manner. Moreover, the stabilizing effect of the feedback can also be used to reduce the external disturbance and minimize the error in such hybrid systems.}

In this paper, we address the null controllability problem for the heat equation posed on a finite interval from a numerical perspective.  We assume that we act on the system by means of a choice of its boundary inputs.  Given $T>0$ and an initial state $u_0$, the aim is to find out whether we can approximate a control (a suitable choice of boundary data) driving the solution from the state $u_0$ to the zero state over the time interval $[0,T]$.  More precisely, we consider the initial boundary value problem on the finite interval $\Omega=(0,L)$:
\begin{eqnarray}\label{heat-fi}
	\begin{cases}
		u_t =  u_{xx}, &x\in \Omega, \ t\in (0,T),\\
		u_x(0,t)=g(t), u_x(L,t)=h(t),  &t\in (0,T),\\
		u(x,0)=u_0(x), &x\in \Omega.
	\end{cases}
\end{eqnarray}
The goal is to compute boundary controls $g(t), h(t)$ so that the solution is steered from the given initial datum $u_0(x)$ to $u(x,T)\equiv 0$ at $t=T$. It is well known that null controllability of this problem can be achieved simply from the right boundary $x=L$; therefore without loss of generality, we will take $g(t)\equiv 0$ and only compute a single control $h(t)$.

{
	The boundary control problem can be placed within a rigorous framework. Indeed, let us consider \eqref{heat-fi} with $T>0$, $u_0\in L^2(\Omega)$, and $h\in L^2(0,T)$. Then, by using the classical tools applied to the wellposedness of nonhomogeneous boundary value problems (see e.g., \cite{Lio72}), one can prove that data to solution map $(u_0,h)\in L^2(\Omega)\times L^2(0,T)\mapsto u\in C([0,T];L^2(\Omega))$ is continuous and therefore $u(T)=u(\cdot,T)$ makes sense as an element of $L^2(\Omega)$ for each given $h\in L^2(0,T)$. Now, define $$R(T,u_0)\equiv \{u(T); u \text{ solves } \eqref{heat-fi} \text{ for some } h\in L^2(0,T)\}.$$ The above set is referred to as the set of reachable states at time $t=T$ starting from the initial state $u_0$ with the aid of $L^2(0,T)$ type boundary controls. We say that \eqref{heat-fi} is \emph{exactly controllable} at time $T$ provided for every initial datum $u_0$, $R(T,u_0)$ coincides exactly with $L^2(\Omega)$.  If $R(T,u_0)$ is only dense in $L^2(\Omega)$, we instead say \eqref{heat-fi} is \emph{approximately controllable}.  Finally, if $0\in R(T,u_0)$, we say that \eqref{heat-fi} is \emph{null-controllable}. Similar terminology is used also for other PDEs and other types of boundary conditions. For time reversible systems such as the wave and Schrödinger equations, it can be shown that exact controllability is equivalent to null-controllability (see e.g., \cite{Zua07}). However, one should not expect a parabolic system such as the heat equation to have the same property. Indeed, heat fails to be exactly controllable due to the smoothing effect. Namely, reachable states are smooth and it is unlikely to steer a system to any desired rough state.  It should be mentioned that controllability properties of a given PDE may change depending on which function spaces are used for states and controls.  Failure of exact controllability can sometimes be fixed with a choice of more regular state space and relatively rougher space for controls. On the other hand, it is well known that null-controllability holds true for the heat equation with respect to $L^2$ framework above on a finite interval (more generally on bounded regular domains in higher dimensions). This fact can be proved through the well-known Hilbert Uniqueness Method \cite{Lio88}.  Moreover, it can be showed that null-controllability for the linear heat equation implies its approximate controllability because eigenfunctions of the Laplacian are reachable from any initial state. One can state the following theorem to summarize the theoretical nature of the control model considered here.}
{
	\begin{thm}\label{controlthm}Let $T>0$ and $u_0\in L^2(\Omega)$. Then, \\
			(i) There exists $h\in L^2(0,T)$ such that the corresponding solution $u$ of \eqref{heat-fi} satisfies $u(\cdot,T)\equiv 0$;\\
			(ii) Given $u_T\in L^2(\Omega)$ and $\epsilon>0$, there exists $h\in L^2(0,T)$ such that the corresponding solution $u$ of \eqref{heat-fi} satisfies $\|u(\cdot,T)-u_T\|_{L^2(\Omega)}<\epsilon$.
	\end{thm}
	\begin{remark}
				Once the solution is steered to the zero state at time $t=T$ (namely once it is stabilized in finite time), then it can be kept at zero state by extending the control as zero for $t>T$. It should be pointed out that the nature of spatial domain is important for controllability. For instance, one can prove for the heat equation that there are initial data which cannot be steered to zero in finite time if the finite domain $\Omega=(0,L)$ is replaced with an unbounded one such as $\Omega=(0,\infty)$ and say with a control acting at $x=0$ (see e.g., \cite{KO20}, \cite{Mi01}).
	\end{remark}
}

An affirmative answer to the control problem is obviously very important due to numerous practical applications.  There is a vast amount of research on the theoretical and numerical approaches for such control problems for partial differential equations (PDEs). Regarding the theoretical approaches see for instance the classical review articles of Russell \cite{Rus78} and Lions \cite{Lio88} and also Zuazua's chapter \cite{Zua07} and the references therein.  For the numerical approaches, we refer the reader to Zuazua's survey article \cite{Zua05}.

Over the last few decades, it became apparent that there is no unified theory regarding the controllability problems for all classes of PDEs. The treatment and results depend on whether one has a time reversible (wave, Schrödinger) vs. time irreversible (heat) equation or a linear vs. nonlinear version of a given type of PDE.

Numerical approximation of the control, which is the main focus of this paper, is a delicate problem. For the wave equation, it is shown in \cite{Zua05} that because of the existence of high-frequency spurious (nonphysical numerical) waves, the control computed for the model which is discretized by a finite difference method does not always approximate the control for the corresponding continuous dynamics well.  Such numerical pathologies diminish for PDEs that have dissipative or dispersive nature, see \cite[Chapter 8]{Zua05}. In particular, for the heat equation due to the high-frequency damping effect it is known that the numerically computed controls will converge.  However, developing those numerical controls is a challenge because of the regularizing effect of the heat kernel. This pathology was studied for the internal control problem in \cite{Mun10} and authors remarked that the situation may be even worse for the boundary control problem that we study here.

Traditionally, controllability problems for PDEs were studied through observability properties of corresponding adjoint systems and authors used some very technical tools such as Carleman estimates, see for instance \cite{Fur96}, \cite{Leb95}.  There are not as much effort given regarding the numerical construction of controls for PDEs and we are aware of only a few works in this direction.  For instance, Rosier \cite{Ros97n} studied the numerical boundary controllability of the Korteweg de-Vries equation on a finite interval by implementing the Hilbert Uniqueness Method (HUM) numerically.   The advantage of HUM is that it allows one to obtain a control which minimizes the $L^2(0,T)$ norm among all controls which steer solutions to desired final state.  Therefore, this approach is desirable if one wants to minimize the cost of control.  In addition, it was numerically shown in \cite{Ros97n} that the $L^2(0,T)$ norm of the control is a nonincreasing function of $T$, suggesting that one does not need to work with large $T$ to achieve optimal controllability.  However, HUM requires one to solve more than one initial boundary value problems numerically, thereby doubling the effect of numerical errors. Namely, one first solves a backward problem and then a forward problem using boundary traces of the backward problem. Moreover,  \cite{Ros97n} shows that the numerical error in HUM gets larger as the frequency of oscillations and points of discontinuities in the target state increase.

On the other hand, direct methods for controllability problems are rather rare. Some earlier theoretical works in this direction for parabolic equations are \cite{Jon77}, \cite{Lar00}, \cite{Guo95}, and \cite{Lit78}. Regarding the direct approaches for constructing numerical approximations of null controls, we are aware of the so called `flatness' approach \cite{Mar14}.  This method relies on constructing the solution and the control at the same time via the derivatives of a flat output. For the problem considered here, this is simply the Dirichlet trace of the solution of \eqref{heat-fi} at $x=0$.  The idea is to utilize the fact that there is a {one-to-one} correspondence between arbitrary Dirichlet traces $y(t)\equiv u(0,t)$ and solutions of \eqref{heat-fi}. One can then prove that the solution and the control can be constructed for \eqref{heat-fi} simultaneously, assuming the solution is of the form $u(x,t)=\sum_{i=0}^\infty\frac{a_i(t)}{i!}x^i$ and the derivatives $y^{(i)}$ are known. In particular, one can deduce $a_{2i+1}=0$, $a_{2i}=\frac{y^{(i)}}{(2i)!}$ and $h(t)=\sum_{i=1}^{\infty}\frac{y^{(i)}}{(2i-1)!}$ if $y^{(i)}$ are known and $u$ is of the given form.  One can ensure that the latter two conditions hold by starting with a suitably chosen and explicit function $y(t)$ of certain Gevrey class.  However, this would imply that the solution itself is also of certain Gevrey class at all times including at the beginning. Namely, this would require initial datum to be also smooth but it is not the case in general. Nevertheless, one can still make the method work. Indeed, one can first steer the solution from the possibly rough initial state $u_0$ to a smoother state by letting the system run with zero boundary input $h(t)=0$ for a short while -say on a time interval $[0,\tau]$ and then supply $h(t)=\sum_{i=1}^{\infty}\frac{y^{(i)}}{(2i-1)!}$ on the interval $(\tau, T]$. The flatness approach is a powerful modern technique in the sense that it is independent of the energy method but its implementation requires some computational effort.  For instance, one needs to be able to numerically compute a large number of derivatives of some functions used in the algorithm. In addition, it is still not known very well how some parameters such as $\tau$ must be chosen for optimal results in the flatness approach.

In this work we propose an alternative approach based on the unified transform, also known as the Fokas
method, in order to determine the boundary control $h(t)$ for the problem
\eqref{heat-fi}. This method was introduced by Fokas in \cite{F97}, for the
analysis of initial and boundary value problems (IBVPs) of integrable nonlinear
PDEs. Later, it led to the emergence of a
novel and powerful approach for studying, both analytically and
numerically, boundary value problems for linear PDEs, \cite{F02, F08}. The
Fokas method consists of two basic elements, (a) the \textit{global
	relation} of the initial condition with the known and the unknown boundary values, and
(b) the \textit{integral representation} of the solution, which contains
both given data and unknown boundary values. A typical procedure for implementing
this method is to use the global relation in order to compute the
contribution of the unknown values in the integral representation,
hence obtain a form of the solution which will contain only the known
initial and boundary conditions of the problem. However, there are many
works that make use of only one of these two elements independently: For
example the power of the global relation has been apparent in problems of
fluid dynamics \cite{AFM06,AF11,FN12,FK17}; also the usage of the integral
representation was crucial for problems related with well-posedness, see e.g.,
\cite{FHM17, OY19}.

In \cite{KO20}, two of the authors implemented the global relation in
order to prove the lack of the null-controllability of the heat equation on
the half line; this is a well known result, but the Fokas
method allowed for an elementary and rather short proof, which was directly
extended to arbitrarily many spatial dimensions. Therein, a
characterisation of the boundary control problem was introduced for the
case of the Dirichlet problem on the finite interval. Here, we derive a
similar characterization stated in \eqref{int-eq-def} for the Neumann problem. Furthermore, we present in detail a numerical scheme which yields
the control $h(t)$; the steps of the numerical scheme are presented in
section \ref{subsecnum}. We performed an extensive study of the performance of the
algorithm, determining also the optimal choice of the relevant parameters;
we summarise these results in section \ref{subsecopt}. We emphasize that for the
implementation of the algorithm we need to discretize the initial condition
$u_0(x)$ into $N+1$ points and compute $(N+1)\times(N+2)$ integrals with
exponentially decaying integrands. Then, inverting an $(N+1)\times(N+1)$
matrix is sufficient for the derivation of the control with an
exponentially small error of order $10^{-2(N-1)}$. The whole procedure for
$N\le10$ takes only few seconds on a laptop.

{ The direct numerical control method presented in this paper relies on the availability of a representation formula obtained through Fokas method, also known as the unified transform method (UTM).  It is well known that UTM can successfully solve a wide range of linear constant coefficient PDEs.  UTM is also applicable to linear systems of equations or to higher-order evolution equations which can be transformed into such systems.  One of the authors, have utilized UTM in controllability of the linear Schr\"odinger equation, \cite{O22}.
	Furthermore, hyperbolic problems (e.g., wave equation, Klein-Gordon equation) can be successfully treated via UTM, see for instance \cite{DGSV18,FK22} and the references therein. There are recent work on the UTM based analysis of Schrödinger and higher order parabolic type equations (\cite{Ba20}, \cite{FHM17}, \cite{BT22}, \cite{OY19}, \cite{OKK22}). Therefore, one can expect to extend our approach to other PDEs by using the associated representation formulas. Furthermore, one of the strongest features of the UTM  is that it provides a unified methodology for the derivation of the integral representation of the solution of IBVPs with different types of boundary conditions, (e.g. Dirichlet, Neumann, Robin, oblique Robin \cite{Ma13}, etc.), thus the computational algorithmic approach described in the current work can be applied to a wide variety of IBVPs.}

\section{Analytical Formulation}\label{secanalytical}
\begin{figure}[!htbp]
	\centering
	\noindent\includegraphics[width=0.5\linewidth]{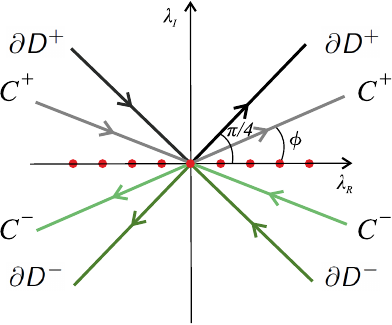}
	\caption{The curves $\partial D^+$ and $\partial D^-$, as well as the curves $C^+$ and $C^-$, with $\phi=\frac{\pi}{8}$.}
	\label{fig1}
\end{figure}

{ In this section we derive the integral equation which characterises the control problem \eqref{heat-fi}.  Theoretical and numerical properties of solutions of PDEs are generally established after one introduces a suitable notion of \emph{solution}. This is crucial especially when the PDE model involves data (initial and/or boundary type) which belong to certain function spaces of low regularity because the classical derivatives in the PDE may not be well-defined or known to exist a priori in the pointwise sense. Therefore, sometimes one needs special and rather weaker formulations to define solutions which do not require smoothness in the classical pointwise sense.  These are generally referred to as \emph{weak} solutions. There are various approches for defining weak solutions. The most classical approach is to interpret the PDE in the distributional sense.  Therefore, one can transfer the derivatives in the PDE onto smooth test functions and look for a solution in the wider space of generalized functions. Another major approach for defining weak solutions is to use \emph{representation formulas}.  Such formulas can be either obtained by using an abstract technique such as semigroup theory or a concrete approach such as an integral transform (Fourier, Laplace, Fokas's UTM, etc). If an integral transform is used, a representation formula is obtained at first assuming the sought after solution is smooth, has nice decay properties and satisfies necessary compatibility conditions.  Once an integral formula is obtained, it can be checked that the same formula is well defined even under much weaker regularity and compatibility properties. This allows one to take the representation formula to be the definition of weak solutions.  {The fact that UTM formula can define low regularity solutions which do not necessarily comply with compatibility was justified in recent articles for certain PDEs, see for instance \cite{OY19}, \cite{HM20}, and \cite{BT22}.} In particular, the integral representation obtained through Fokas method in this paper follows the same idea. Fokas's method \cite{F08, FKinpress, DTV14}  yields the integral representation of the solution of \eqref{heat-fi} given by:}
\begin{align}\label{sol-1}
	u&(x,t)=
	\int_{-\infty }^{\infty }e^{i\lambda x-\lambda^2 t}\hat{u}_{0}(\lambda )\frac{d\lambda}{2\pi}
	\notag -\int_{\partial D^{+}} \frac{e^{i\lambda x-\lambda^2 t}}{e^{i\lambda L} - e^{-i\lambda L}}
	\left[ e^{i\lambda L} \hat{u}_{0}(\lambda ) + e^{-i\lambda L} \hat{u}_{0}(-\lambda )\right.\notag\\
	&\left.-2e^{-i\lambda L}\tilde{g}(\lambda ^{2},t) +2\tilde{h}(\lambda ^{2},t)
	\right] \frac{d\lambda}{2\pi}
	\notag -\int_{\partial D^{-}} \frac{e^{i\lambda x-\lambda^2 t}}{e^{i\lambda L} - e^{-i\lambda L}}
	\left[ e^{-i\lambda L} \hat{u}_{0}(\lambda ) + e^{-i\lambda L} \hat{u}_{0}(-\lambda )\right.\notag\\
	&\left.-2e^{-i\lambda L}\tilde{g}(\lambda ^{2},t) +2\tilde{h}(\lambda ^{2},t)
	\right] \frac{d\lambda}{2\pi} , \ \ x\in[0,L], \ t>0,
\end{align}
where $\partial D^{+}$ and $\partial D^{-}$ are depicted in Figure \ref{fig1}, and for $\lambda\in\C$,
$$\hat{u}_0(\lambda)=\int_0^L e^{-i \lambda x} u_0(x)dx, \qquad \tilde{g}(\lambda ^{2},t)=\int_0^t e^{\lambda^2 s} g(s)ds,$$
$$\tilde{h}(\lambda ^{2},t)=\int_0^t e^{\lambda^2 s} h(s)ds. $$

{\begin{remark}
		The curve $\partial D^+$ is formed by the two rays $\left( e^\frac{3i\pi}{4} \infty,0\right)\cup\left(0,e^\frac{i\pi}{4} \infty\right)$. However, since $\lambda=0$ may be a singular point, for matters of rigour in what follows we equip this curve with a quarter-circle of small radius $\epsilon>0$, around $\lambda=0$. This deformation is doable in \eqref{IR1}, without loss of contribution as $\epsilon\to0$. In analogy we define $C^+$, which is formed by the two rays $\left( e^\frac{7i\pi}{8} \infty,0\right)\cup\left(0,e^\frac{i\pi}{8} \infty\right)$, equipped with 3/8 of a circle of small radius, around $\lambda=0$. The curves $\partial D^-$ and $C^-$ are defined as the opposite of  $\partial D^+$ and $C^+$, respectively.
\end{remark}}

For matters of presentation we sketch the derivation of \eqref{sol-1}. We follow the 3 steps of the general methodology of the Fokas method, as it is described for example in Chapter 10 of \cite{FKinpress}. The first step involves the derivation of the Global Relation (Eq. (10.1) of \cite{FKinpress}) , namely
\begin{align} \label{GR1}
	\hspace*{-7mm} e^{\lambda ^{2}t}\hat{u}(\lambda ,t)=\hat{u}_{0}(\lambda )
	-\tilde{g}_{1}(\lambda ^{2},t)-i\lambda\tilde{g}_{0}(\lambda ^{2},t)
	 +e^{-i\lambda L}\left[
	\tilde{h}_{1}(\lambda ^{2},t)+i\lambda\tilde{h}_{0}(\lambda ^{2},t)
	\right] , \quad \lambda \in\mathbb{C},
\end{align}
where $\hat{u}$ and $\hat{u}_{0}$ are the finite Fourier transforms of $u(x,t)$ and $u_{0}(x)$, defined by
\begin{align*}
	\hat{u}(\lambda ,t)=\int_{0}^{L}e^{-i\lambda x}u(x,t)dx, \quad \hat{u}_{0}(\lambda )=\int_{0}^{L}e^{-i\lambda x}u_{0}(x)dx, \qquad \lambda \in\mathbb{C},
\end{align*}
and
\begin{equation*}
	\tilde{f}(\lambda ,t)=\int_{0}^{t}e^{\lambda s }f(s )ds , \qquad t>0,\  \lambda \in\mathbb{C},
\end{equation*}
with $g_0(t )=u(0,t)$, $h_0(t )=u(L,t),  \  t>0$.\\
\indent The second step involves the derivation of the Integral Representation of the solution  (Eq. (10.3a) of \cite{FKinpress}), namely
\begin{align}\label{IR1}
	u(x,t)=&\frac{1}{2\pi}\int_{-\infty }^{\infty }e^{i\lambda x-\lambda ^{2}t}\hat{u}_{0}(\lambda )d\lambda
	-\frac{1}{2\pi}\int_{\partial D^{+}}e^{i\lambda x-\lambda ^{2}t}
	\left[ \tilde{g}_{1}+i\lambda \tilde{g}_{0}   \right] d\lambda
	\notag\\
	&-\frac{1}{2\pi}\int_{\partial D^{-}}e^{-i\lambda (L-x)-\lambda ^{2}t}
	\left[ \tilde{h}_{1}+i\lambda \tilde{h}_{0}   \right] d\lambda.
\end{align}
This equation is obtained by applying the inverse Fourier transform formula to \eqref{GR1}, and  deforming from $\R$ to $\partial D^{+}$ in the integral involving $\tilde{g}_j$, and from $\R$ to $\partial D^{-}$ in the integral involving $\tilde{h}_j$. { This deformation is possible because (a) the integrand of the second integral is bounded and analytic in the domain defined by $\R$ and $\partial D^{+}$; (b) the integrand of the third integral is bounded and analytic in the domain defined by $\R$ and $\partial D^{-}$; (c) both integrands are $O\left(\frac{1}{\lambda}\right), \ \lambda\to\infty$. Indeed,  by using Cauchy's theorem, enhanced with the Jordan's lemma, the above argumentation implies that the integration of the second integrand on the boundary of the domain defined by $\R$ and $\partial D^{+}$ yields zero contribution. Similarly, the integration of the third integrand on the boundary of the domain  defined by $\R$ and $\partial D^{-}$ yields zero contribution.}
\\
\indent In the third step we create one more Global Relation by the transform $\lambda\to -\lambda$ in \eqref{GR1}, which leaves  $\tilde{g}_j$ and  $\tilde{h}_j$ invariant. { Then, we solve these two Global Relations for $\tilde{g}_0$ and  $\tilde{h}_0$ and substitute the resulting expressions in \eqref{IR1}. We note that employing Cauchy's theorem and Jordan's lemma, the integrals involving $\hat{u}(\pm\lambda ,t)$ along $\partial D^\pm$ vanish; denoting $g_1=g$ and $h_1=h$ yields \eqref{sol-1}.}\\
\indent Let $g\equiv 0, \ u(x,T)=0$.
Then, \eqref{sol-1} becomes
\begin{align}\label{sol-2}
	&0=u(x,T)=
	\int_{-\infty }^{\infty }e^{i\lambda x-\lambda^2 T}\hat{u}_{0}(\lambda )\frac{d\lambda}{2\pi}
	\notag \\
	&-\int_{\partial D^{+}} \frac{e^{i\lambda x-\lambda^2 T}}{e^{i\lambda L} - e^{-i\lambda L}}
	\left[ e^{i\lambda L} \hat{u}_{0}(\lambda ) + e^{-i\lambda L} \hat{u}_{0}(-\lambda )\right.\notag\left.+2\tilde{h}(\lambda ^{2},T)
	\right] \frac{d\lambda}{2\pi}
	\notag \\
	&-\int_{\partial D^{-}} \frac{e^{i\lambda x-\lambda^2 T}}{e^{i\lambda L} - e^{-i\lambda L}}
	\left[ e^{-i\lambda L} \hat{u}_{0}(\lambda ) + e^{-i\lambda L} \hat{u}_{0}(-\lambda )\right.\left.+2\tilde{h}(\lambda ^{2},T)
	\right] \frac{d\lambda}{2\pi} , \qquad x\in[0,L].
\end{align}

We first rewrite \eqref{sol-2} in the form
\begin{align*}
	&2\int_{\partial D^{+}} \frac{e^{i\lambda x-\lambda^2 T}}{e^{i\lambda L} - e^{-i\lambda L}} \tilde{h}(\lambda ^{2},T) \frac{d\lambda}{2\pi}\notag+ 2\int_{\partial D^{-}} \frac{e^{i\lambda x-\lambda^2 T}}{e^{i\lambda L} - e^{-i\lambda L}} \tilde{h}(\lambda ^{2},T) \frac{d\lambda}{2\pi}
	\notag\\
	&=\int_{-\infty }^{\infty }e^{i\lambda x-\lambda^2 T}\hat{u}_{0}(\lambda )\frac{d\lambda}{2\pi}
	\notag -\int_{\partial D^{+}} \frac{e^{i\lambda x-\lambda^2 T}}{e^{i\lambda L} - e^{-i\lambda L}}
	\left[ e^{i\lambda L} \hat{u}_{0}(\lambda ) + e^{-i\lambda L} \hat{u}_{0}(-\lambda )
	\right] \frac{d\lambda}{2\pi}
	\notag \\
	&-\int_{\partial D^{-}} \frac{e^{i\lambda x-\lambda^2 T}}{e^{i\lambda L} - e^{-i\lambda L}}
	\left[ e^{-i\lambda L} \hat{u}_{0}(\lambda ) + e^{-i\lambda L} \hat{u}_{0}(-\lambda )
	\right] \frac{d\lambda}{2\pi} .
\end{align*}
Letting $\lambda\to-\lambda$ in the second integral of the left hand side (LHS), which transforms $\partial D^- \to \partial D^+$, we obtain the expression
\begin{align*}
	&2\int_{\partial D^{+}} \frac{e^{i\lambda x-\lambda^2 T}}{e^{i\lambda L} - e^{-i\lambda L}} \tilde{h}(\lambda ^{2},T) \frac{d\lambda}{2\pi}+ 2\int_{\partial D^{-}} \frac{e^{i\lambda x-\lambda^2 T}}{e^{i\lambda L} - e^{-i\lambda L}} \tilde{h}(\lambda ^{2},T) \frac{d\lambda}{2\pi} \\
	&=2\int_{\partial D^{+}} \frac{e^{i\lambda x}+e^{-i\lambda x}}{e^{i\lambda L} - e^{-i\lambda L}}e^{-\lambda^2 T} \tilde{h}(\lambda ^{2},T) \frac{d\lambda}{2\pi}\\ &=-\frac{1}{\pi}\int_{\partial D^+} \cos(\lambda x)\frac{ i e^{-\lambda^2 T} \tilde{h}(\lambda ^{2},T)}{\sin(\lambda L) } d\lambda.
\end{align*}

Also we will deform $\partial D^+$ to $C^+$, which is forming angles of $\frac{\pi}{8}$ with the real axis, in contrast to $\frac{\pi}{4}$ of the rays $\partial D^+$. The curves $C^+$ and $C^-$ are depicted in Figure \ref{fig1}. This deformation is possible because the integrand is bounded and analytic in the region where $\text{Im} \lambda >0$ and $\text{Re} \lambda^2>0$.   { The choice of the curves $C^\pm$ is led by the fact that
	$\phi=\frac{\pi}{8}$ yields the largest growth in the value of the expression $\text{Re} \lambda^2$, hence the integrand has its strongest exponential decay, as $|\lambda|\rightarrow \infty$. This observation will be useful in the computational part of the current work, where the numerical computation of the latter integrals is much quicker and robust than the computation of integrals with oscillatory integrands.\\
	\indent { We note that in general the singularity $\lambda=0$ yields the extra constant contribution $\frac{1}{2L}\tilde{h}(0,T)$; this contribution can be neglected if $\frac{1}{2L}\tilde{h}(0,T)=0$, hence $\lambda=0$ becomes a removable singularity. In what follows for reasons of simplicity of presentation we restrict our computation to control $h(t)$ that satisfy the condition $\tilde{h}(0,T)=\int_0^T h(t)dt=0$.}
	
	
	Thus, equation \eqref{sol-2} takes the form
	\begin{align}\label{int-eq-def}
		-\int_{C^+} \cos(\lambda x)\frac{ i e^{-\lambda^2 T} \tilde{h}(\lambda ^{2},T)}{\sin(\lambda L)} d\lambda
		=P(x;L,T), \quad x\in[0,L],
	\end{align}
	where
	\begin{align}\label{def:P}
		&P(x;L,T)=\frac{1}{2}\int_{-\infty }^{\infty }e^{i\lambda x-\lambda^2 T}\hat{u}_{0}(\lambda )d\lambda
		\notag \\
		&-\frac{1}{2}\int_{\partial D^{+}} \frac{e^{i\lambda x-\lambda^2 T}}{e^{i\lambda L} - e^{-i\lambda L}}
		\left[ e^{i\lambda L} \hat{u}_{0}(\lambda ) + e^{-i\lambda L} \hat{u}_{0}(-\lambda )
		\right] d\lambda
		\notag \\
		&-\frac{1}{2}\int_{\partial D^{-}} \frac{e^{i\lambda x-\lambda^2 T}}{e^{i\lambda L} - e^{-i\lambda L}}
		\left[ e^{-i\lambda L} \hat{u}_{0}(\lambda ) + e^{-i\lambda L} \hat{u}_{0}(-\lambda )
		\right] d\lambda .
	\end{align}
	
	The control $h(t)$ is obtained as a solution of the integral equation \eqref{int-eq-def}. We note that the right hand side (RHS) of \eqref{int-eq-def} is known, and that we need to exploit the fact that this equation is valid for \textit{all} $x\in[0,L]$.

	\section{Construction of the numerical algorithm}\label{secconstruct}
	
	{ In this section we prepare the apparatus for computing numerically the control $h(t)$. We emphasise on the fact that the exponential decay of the integrands involved in the numerical computation allows for almost arbitrary choice of the basis  $\displaystyle\{\phi_n(t)\}_{n=1}^{\infty}$, which we use to approximate the control $h(t)$ below. For this reason, in what follows we chose one of the most well-known basis, namely the sine-Fourier with the only restriction to vanish at $t=\tau$ and $t=T$.}
	
	\subsection{Rewriting \eqref{int-eq-def} via a basis for $h(t)$}\label{subsecrew}
	
	Let the basis $\displaystyle\{\phi_n(t)\}_{n=1}^{\infty}$, and the approximation
	\begin{equation} \label{def:h-basis}
		h(t)=\sum_{n=1}^{N+1} a_n \phi_n(t).
	\end{equation}
	Then,
	\begin{equation} \label{def:ht-basis}
		\tilde{h}(\lambda^2,t)=\sum_{n=1}^{N+1} a_n \int_0^t e^{\lambda^2 s} \phi_n(s)ds.
	\end{equation}
	Thus,  by fixing $L$ and $T$, the integral equation \eqref{int-eq-def} takes the form
	\begin{align}\label{int-eq-g-num-2}
		\sum_{n=1}^{N+1} a_n F_n(x)=G(x) .
	\end{align}
	The definitions of the above functions are given as follows.
	\begin{itemize}
		\item Definition of $F_n$:
		\begin{equation}\label{def-Fn}
			F_n(x)=-\int_{C^+} \cos(\lambda x)\frac{i e^{-\lambda^2 T} B_n(\lambda)}{\sin(\lambda L) } d\lambda
		\end{equation}
		with
		\begin{equation}\label{def-Bn}
			B_n(\lambda)= \int_0^{T} e^{\lambda^2 s} \phi_n(s)ds.
		\end{equation}
		\item Definition of $G$:
		\begin{equation}\label{def-G}
			G(x)=P(x;L,T), \qquad x\in[0,L].
		\end{equation}
	\end{itemize}

	In all the examples we have chosen the functions $\big\{\phi_n(t)\big\}_{n=1}^{N+1}$ of the general sine-Fourier basis supported on the interval $t\in(\tau,T)$, namely
	\begin{equation}\label{def-phi-num-g}
		\phi_n(t) =\begin{cases} 0, & 0\le t\le \tau,\\
			\sin \left( \pi  n \dfrac{ t-\tau}{T-\tau}\right),  & \tau\le t \le T,\\
			0, & T<t.\end{cases}
	\end{equation}
	Then \eqref{def-Bn} yields
	\begin{align}\label{def-BN-num-g}
		B_n(\lambda)&{=\int_\tau^T e^{\lambda ^2 s} \sin \left(\frac{\pi  n (s-\tau)}{T-\tau}\right) \, ds} \notag\\
		&{=\frac{\lambda ^2 (\tau-T)^2 e^{\lambda ^2 s} \sin \left(\frac{\pi  n (\tau-s)}{\tau-T}\right)}{\lambda ^4 (\tau-T)^2+\pi ^2 n^2} \Bigg|_{s=\tau}^T}
		{ +\frac{\pi  n (\tau-T) e^{\lambda ^2 s} \cos \left(\frac{\pi  n (\tau-s)}{\tau-T}\right)}{\lambda ^4 (\tau-T)^2+\pi ^2 n^2}\Bigg|_{s=\tau}^T} \notag\\
		&=\frac{\pi  n (\tau-T) \left[ (-1)^n e^{\lambda ^2 T}-e^{\tau \lambda ^2}\right]}{\lambda ^4 (\tau-T)^2+\pi ^2 n^2}.
	\end{align}

	\subsection{The steps of the algorithm for solving \eqref{int-eq-g-num-2}}\label{subsecnum}
	\begin{itemize}
		\item[1.] Fix $L$ and $T$.  Choose $\displaystyle\{x_k\}_{k=0}^{K}$ in the interval $x_k\in[0,L]$.
		\item[2.] Compute numerically $F_n(x_k)$ via \eqref{def-Fn} and \eqref{def-BN-num-g}. Then, the LHS of \eqref{int-eq-g-num-2} yields a $(K+1)\times (N+1)$ matrix, which we denote as $\mathcal{F}$, with elements $F_{n,k}=F_n(x_k)$.
		\item[3.] Compute numerically $G(x_k)$ via \eqref{def-G} and \eqref{def:P}.  This constructs a $(K+1)$ vector, which  we denote as $\mathfrak{g}$, with elements $G_{k}=G(x_k)$.
		\item[4.] Solve the equation
		\begin{equation}\label{int-eq-c-num}
			\mathcal{F}\cdot\alpha=\mathfrak{g},
		\end{equation}
		for $\alpha$, which is a $(N+1)$ vector with entries the unknowns $a_n$. Then, \eqref{def:h-basis} yields the control.
	\end{itemize}
	{
		\begin{remark}
					Note that the last step in above algorithm requires $\mathcal{F}^{-1}$ to exist.   This can always be achieved because from definition \eqref{def-Fn} we observe that the entries of $\mathcal{F}$ given by $F_n(x_k)$ can be modified if necessary by choosing sligtly different $x_k$. Namely, if $\mathcal{F}$ is singular for an initial choice of $x_k$'s, we can modify one or more $x_k$ values to make $\mathcal{F}$ nonsingular as the only dependence of  $F_n(x_k)$ on $x_k$ in the integral at the RHS of \eqref{def-Fn} is the cosine term. In practice, we did not have to make such modification for any of the numerical simulations.
		\end{remark}
		\begin{remark}
					The flatness approach of \cite{Mar14} first constructs an exact null control (as an infinite series) and then takes the approximate null control as a finite truncation of this series.  Then, it is shown via norm estimates that the difference between the exact  solution with exact null control and the numerical solution with truncated null control is not large. In contrast, as our algorithm suggests, we propose a direct numerical approximation of the null control function.  The approximate numerical control function that we construct has by definition a truncated nature (a finite series) and we show in subsequent sections that the system behaves in a manner that solution vanishes numerically at time $t=T$ with an exponenentially small error even for a rather small choice of $N$.
		\end{remark}
	}
	
	\section{First Specific Example}\label{secfirst}
	
	{In this section, we evaluate $P(x;L,T)$ via the RHS of \eqref{def:P} for $L=1$,  $T=1/2$, and  the following initial condition:
	}
	
	Let $u_0(x)=\begin{cases} -1, &x\in(0,1/2) \\ 1, &x\in(1/2,1)  \end{cases}$.  Then,
	\begin{equation}\label{hat-u-0}
		\hat{u}_0(\lambda)=\frac{i e^{-i \lambda } \left(-1+e^{\frac{i \lambda }{2}}\right)^2}{\lambda }=\frac{i}{\lambda }-\frac{2 i e^{-\frac{1}{2} (i \lambda )}}{\lambda }+\frac{i e^{-i \lambda }}{\lambda }.
	\end{equation}
	{ It turns out that the RHS of \eqref{int-eq-g-num-2} (namely the expression of $G(x)=P(x;1,T)$, with  $P$ given by \eqref{def:P}) can be simplified to the RHS of \eqref{G-example1}.\\
		\indent Indeed, employing the transform $\lambda\to-\lambda$ in \eqref{hat-u-0} and employing them in the integrands of the second and third integrals of \eqref{def:P} yields
		\begin{align}\label{help-r-1}
			-\frac{1}{2} \frac{e^{i\lambda } \hat{u}_{0}(\lambda ) + e^{-i\lambda } \hat{u}_{0}(-\lambda )}{e^{i\lambda } - e^{-i\lambda }} = -\frac{i}{2\lambda}\left[1-\frac{1}{\cos\left(\frac{\lambda}{2}\right)}\right]
		\end{align}
		and
		\begin{align}\label{help-r-2}
			-\frac{1}{2} \frac{e^{-i\lambda } \hat{u}_{0}(\lambda ) + e^{-i\lambda } \hat{u}_{0}(-\lambda )}{e^{i\lambda } - e^{-i\lambda }} = \frac{i e^{-i\lambda }}{2\lambda}\left[1-\frac{1}{\cos\left(\frac{\lambda}{2}\right)}\right],
		\end{align}
		respectively. Equations \eqref{help-r-1} and \eqref{help-r-2} imply also that $\lambda=0$ is a removable singularity for the last two integrals in \eqref{def:P}, hence we are allowed to deform the  curves $\partial D^\pm$ to $C^\pm$, using Cauchy's theorem without the contribution of any residue contribution.\\
		\indent By employing \eqref{hat-u-0} in the first term of the RHS of \eqref{def:P}, we obtain three integrals. Due to the analyticity and boundedness of the integrands, the first one can be deformed to the $C^+$, whereas the third one  can be deformed to the $C^-$. Furthermore, for $x<1/2$, the second term can be deformed to the $C^-$; for $x>1/2$, the second term can be deformed to the $C^+$. Of course this procedure yields an additional contribution of the residue due to the pole $\lambda=0$, which we calculate.
		The curves $C^+$ and $C^-$ are forming angles of $\phi=\frac{\pi}{8}$ and are depicted in Figure \ref{fig1}. }

	Indeed,
	\begin{itemize}
		\item for $x\in(0,1/2)$ the first term of the RHS of \eqref{def:P} reads 
		\begin{align}\label{help-r-3}
			&\frac{1}{2}\int_{-\infty }^{\infty }e^{i\lambda x-\lambda^2 T}\hat{u}_{0}(\lambda )d\lambda \notag \\
			=&\frac{1}{2}\int_{-\infty }^{\infty } e^{i\lambda x-\lambda^2 T}\left( \frac{i}{\lambda }-\frac{2 i e^{-\frac{1}{2} (i \lambda )}}{\lambda }+\frac{i e^{-i \lambda }}{\lambda } \right) d\lambda \notag \\
			=&\frac{1}{2}\int_{C^+} e^{i\lambda x-\lambda^2 T}\left( \frac{i}{\lambda }\right) d\lambda -\frac{1}{2}\int_{C^-} e^{i\lambda x-\lambda^2 T}\left( -\frac{2 i e^{-\frac{1}{2} (i \lambda )}}{\lambda }+\frac{i e^{-i \lambda }}{\lambda } \right) d\lambda \notag \\
			&+ \frac{\pi}{8} i  \text{Res} \left\{e^{i\lambda x-\lambda^2 T}\left( \frac{i}{\lambda }\right)  \right\}_{\lambda=0} - \frac{\pi}{8} i \text{Res} \left\{e^{i\lambda x-\lambda^2 T}\left( -\frac{2 i e^{-\frac{1}{2} (i \lambda )}}{\lambda }+\frac{i e^{-i \lambda }}{\lambda } \right) \right\}_{\lambda=0} \notag\\
			=&\frac{1}{2}\int_{C^+} e^{i\lambda x-\lambda^2 T}\left( \frac{i}{\lambda }\right) d\lambda+\frac{1}{2}\int_{C^-} e^{i\lambda x-\lambda^2 T}\left( \frac{2 i e^{-\frac{1}{2} (i \lambda )}}{\lambda }-\frac{i e^{-i \lambda }}{\lambda } \right) d\lambda
			- \frac{\pi}{4}.
		\end{align}
		\item similarly for $x\in(1/2,1)$ the first term of the RHS of \eqref{def:P} reads
		\begin{align*}
			&\frac{1}{2}\int_{-\infty }^{\infty }e^{i\lambda x-\lambda^2 T}\hat{u}_{0}(\lambda )d\lambda\\
			=&\frac{1}{2}\int_{C^+} e^{i\lambda x-\lambda^2 T}\left( \frac{i}{\lambda } -\frac{2 i e^{-\frac{1}{2} (i \lambda )}}{\lambda }\right) d\lambda +\frac{1}{2}\int_{C^-} e^{i\lambda x-\lambda^2 T}\left(-\frac{i e^{-i \lambda }}{\lambda } \right) d\lambda + \frac{\pi}{4}.
		\end{align*}
	\end{itemize}
	
	{ Thus, by applying  \eqref{help-r-1}, \eqref{help-r-2} and \eqref{help-r-3} in \eqref{def:P}, } and deforming the curves $\partial D^{\pm}$ to $C^{\pm}$ respectively, straightforward calculations yield the following expression
	\begin{align}\label{int-eq-s}
		P(x;1,T)
		=&\frac{1}{2}\int_{C^+} e^{i\lambda x-\lambda^2 T} \frac{i}{\lambda  \cos \left(\frac{\lambda }{2}\right)} d\lambda+\frac{1}{2}\int_{C^-} e^{i\lambda x-\lambda^2 T} \frac{i}{\lambda  \cos \left(\frac{\lambda }{2}\right)} d\lambda - \frac{\pi}{4}, \notag\\
		=&\int_{C^+}\cos(\lambda x)  \frac{ i e^{-\lambda^2 T}}{\lambda  \cos \left(\frac{\lambda }{2}\right)} d\lambda - \frac{\pi}{4}, \qquad x\in(0,1/2),
	\end{align}
	where the last equality was derived by letting $\lambda\to-\lambda$ in the second integral of the RHS, which transforms $C^- \to  C^+$.
	Similarly,
	\begin{align}\label{int-eq-b}
		P(x;1,T)
		=-\int_{C^+}\cos\big(\lambda (x-1)\big)  \frac{ i e^{-\lambda^2 T} }{\lambda  \cos \left(\frac{\lambda }{2}\right)} d\lambda + \frac{\pi}{4},
		\quad x\in(1/2,1).
	\end{align}

	
	By fixing $T=1/2$, equations \eqref{int-eq-s} and \eqref{int-eq-b} yield
	\begin{equation}\label{G-example1}
		G(x)=\begin{cases}U(x), & x\in(0,1/2) \\
			-U(1-x), & x\in(1/2,1) \end{cases} ,
	\end{equation}
	with $$U(x)= \int_{C^+}\cos(\lambda x)  \frac{ i e^{-\frac{\lambda^2}{2}}}{\lambda
		\cos \left(\frac{\lambda }{2}\right)} d\lambda -  \frac{\pi}{4}.$$
	Thus, the integral equation \eqref{int-eq-def} takes the form
	\begin{align}\label{int-eq-g}
		-\int_{C^+} \cos(\lambda x)\frac{i e^{-\frac{\lambda^2}{2}} \tilde{h}\left(\lambda ^{2},\frac{1}{2}\right)}{\sin\lambda } d\lambda
		=G(x), \ \ x\in[0,1],
	\end{align}
	with $\tilde{h}$ defined by \eqref{def:ht-basis}.

	\section{Implementation of the algorithm}\label{secimp}
	
	{ This section demonstrates the performance of our algorithm for different values of the main parameters of the problem. Our goal is to show the robustness of this implementation, which yields qualitatively the same results for \textit{all} our simulations.
	}
	
	\subsection{Writing \eqref{int-eq-g} in the form  of \eqref{int-eq-c-num}}\label{sebsecwriting}
	
	Let $N=K$, and $x_k=\frac{k}{N}, \ k=0,\ldots,N$.

	\subsubsection*{Evaluate $\mathfrak{g}$}\label{sebseceval}
	
	We evaluate $\mathfrak{g}$ by employing \eqref{G-example1}, namely
	$$G_k=G(x_k)=G\left(\frac{k}{N}\right), \qquad k=0,\ldots,N.$$
	By employing the numerical integration
	\begin{align} \label{def-U-num-2}
		U(x)= \int_0^{+\infty} \text{Re} \left[ \cos( r e^{i\pi/8} x)  \frac{i  e^{-\frac{ r^2 e^{i\pi/4}}{2}}}{ r
			\cos \left(\frac{ r e^{i\pi/8} }{2}\right)}\right. \left.   -\cos( r e^{7i\pi/8} x)  \frac{i e^{-\frac{ r^2 e^{7i\pi/4}}{2}}}{ r
			\cos \left(\frac{ r e^{7i\pi/8} }{2}\right)} \right]  dr  - \frac{\pi}{4},
	\end{align}
	we obtain
	\begin{align}\label{def-G-num-1}
		G_k= \begin{cases} U\left(\frac{k}{N}\right), & k=0,\ldots, \left[\frac{N}{2}\right],\\
			- U\left(1-\frac{k}{N}\right), & k= \left[\frac{N}{2}\right]+1,\ldots,N.
		\end{cases}
	\end{align}
	
	
	\subsubsection*{Evaluate $\mathcal{F}$}\label{sebseceval2}
	
	Choose the basis $\left\{\phi_n(s)\right\}_{n=1}^{N+1}$ via \eqref{def-phi-num-g}. Then $\left\{B_n(\lambda)\right\}_{n=1}^N$ is given by \eqref{def-BN-num-g}.
	
	Then, for $n=1,\ldots,N+1 \ $ and $ \ k=0,\ldots,N, \ $ the entries  $F_{n,k}$ are numerically evaluated by
	$F_{n,k}$ are numerically evaluated by
	\begin{align}\label{def-Fn-num-3}
		&F_{n,k}=F_n\left(\frac{k}{N}\right)\notag\\
		=&-\int_0^{+\infty}  \text{Re} \Bigg[ \cos\left( r e^{i\pi/8} \frac{k}{N}\right)\frac{ i e^{-\frac{ r^2 e^{i\pi/4}}{2}} B_n( r e^{i\pi/8})}{\sin (r e^{i\pi/8}) }   e^{i\pi/8} \notag \\
		&- \cos\left( r e^{7i\pi/8} \frac{k}{N}\right)\frac{ i e^{-\frac{ r^2 e^{7i\pi/4}}{2}} B_n( r e^{7i\pi/8})}{\sin (r e^{7i\pi/8}) }   e^{7i\pi/8} \Bigg] dr.
	\end{align}
	
	\subsubsection*{Solving \eqref{int-eq-c-num}, namely $\mathcal{F}\cdot \alpha =\mathfrak{g}$}\label{subsecsolv}
	
	We first present the particular example $\tau=0.3$, which is identical with the one introduced in \cite{Mar13}.
	
	For $N=8$ we obtain
	\begin{equation}\label{an-ex}
		\left(\begin{array}{c} a_1 \\ a_2 \\ a_3 \\ a_4 \\ a_5 \\ a_6 \\ a_7 \\ a_8 \\a_9\\
		\end{array}\right)
		=\alpha=\mathcal{F}^{-1}\cdot \mathfrak{g}=
		\left(\begin{array}{c} -0.43685 \\ -0.72935 \\ -0.42262 \\ 0.69991 \\ 1.9004 \\ 2.1164 \\ 1.3298 \\ 0.46097 \\ 0.068970\\
		\end{array}\right)
	\end{equation}
	
	Then,
	\begin{equation}\label{h1-sol-ex}
		h(t)=\sum_{n=1}^9 a_n \phi_n(t),
	\end{equation}
	with $a_n$ given in \eqref{an-ex} and $\phi_n(t)$ given in \eqref{def-phi-num-g} is depicted in Figure \ref{fig3} by the red curve. Also the $L^2$-norm, defined by $\left(\int_0^t h^2(s)ds\right)^{1/2}$ is depicted in Figure \ref{fig3} by the blue curve.
	\begin{figure}[!htbp]
		\centering
		\noindent\includegraphics[width=0.6\linewidth ]{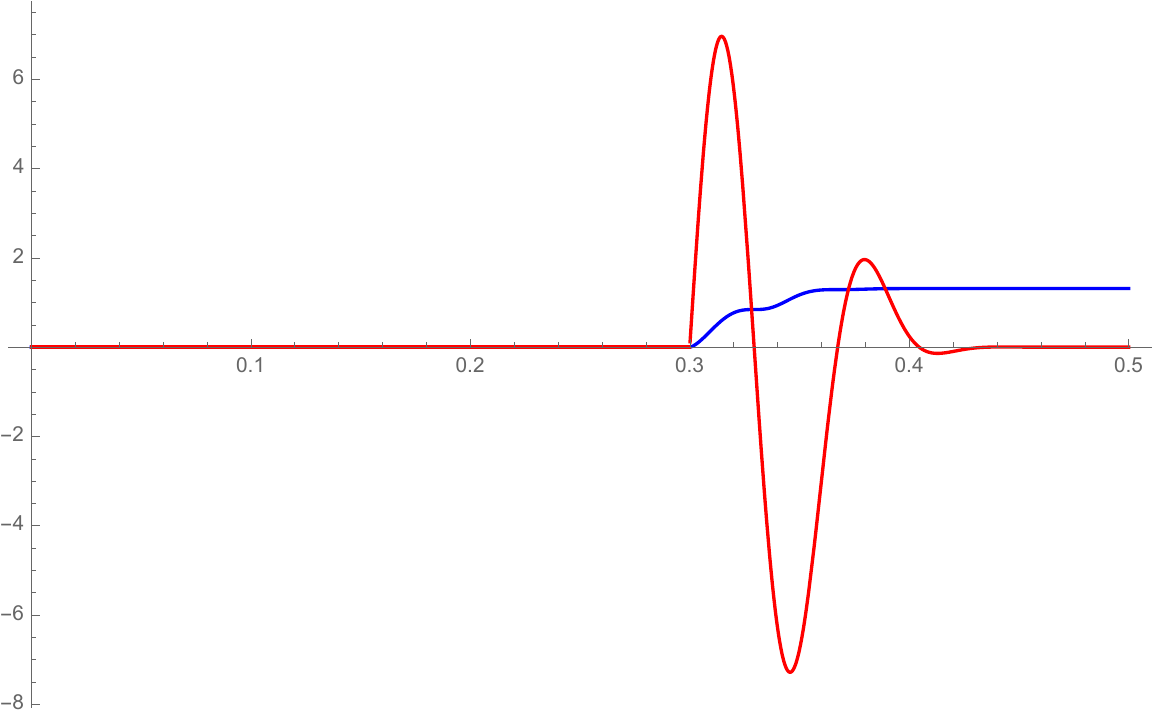}
		\caption{The control  $h(t)$ (red) and its $L^2$-norm (blue).}
		\label{fig3}
	\end{figure}
	
	\subsubsection*{Verification}\label{subsecver}
	
	For matters of completeness we evaluate $u(x,t)$ via \eqref{sol-1}, by making the proper substitutions. Namely, employ the solution
	\begin{equation}\label{sol-ver-0}
		u(x,t)= \frac{1}{\pi}P(x;1,t) + \frac{1}{\pi}\int_{C^+} \cos(\lambda x)\frac{i e^{-\lambda^2 t} \tilde{h}\left(\lambda ^{2},t\right)}{\sin\lambda } d\lambda,
	\end{equation}
	where
	\begin{equation}\label{P-example1-v}
		P(x;1,t)=\begin{cases}Q(x,t), & x\in(0,1/2) \\
			-Q(1-x,t), & x\in(1/2,1) \end{cases} ,
	\end{equation}
	with
	\begin{align}\label{Q-example1-v}
		Q(x,t)
		=\int_{C^+}\cos\big(\lambda x\big)  \frac{ i e^{-\lambda^2 t} }{\lambda  \cos \left(\frac{\lambda }{2}\right)} d\lambda - \frac{\pi}{4}, \qquad x\in(0,1/2).
	\end{align}
	Furthermore, $\ \displaystyle \tilde{h}\left(\lambda ^{2},t\right)=\int_0^t e^{\lambda^2 s} h(s)ds\ $ and $h(s)$ given by \eqref{h1-sol-ex}.
	Namely,
	$$ \tilde{h}\left(\lambda ^{2},t\right)=\sum_{n=1}^{N+1}a_n b_n(\lambda,t),$$
	with $b_n(\lambda,t)$ defined as follows
	\begin{equation}\label{def-bnt-num-g}
		b_n(\lambda,t)= \begin{cases} 0, & 0<t<\tau, \\
			\dfrac{(T-\tau) }{\lambda ^4 (\tau-T)^2+\pi ^2 n^2} \left\{\pi  n e^{\lambda ^2 \tau}\right. \\ \left. -e^{\lambda ^2 t} \left[\lambda ^2 (\tau-T) \sin \left(\pi  n \frac{ t-\tau}{T-\tau}\right) \right. \right. \\ \left. \left. +\pi  n \cos \left( \pi  n \frac{t-\tau}{T-\tau} \right)\right]\right\},
			& \tau\le t \le T\\
			0, & t>T. \end{cases}
	\end{equation}
	
	Obviously, $b_n(\lambda,T)=B_n(\lambda)$.
	
	This solution is depicted in Figure \ref{fig4}.  Furthermore $||u(x,T)||_{L^2}\approx 2.32 \times 10^{-11}$, namely $u(x,T)$ vanishes up to an error of order $10^{-11}$.

	\begin{figure}[!htbp]
		
		\centering
		\noindent\includegraphics[width=0.6\linewidth ]{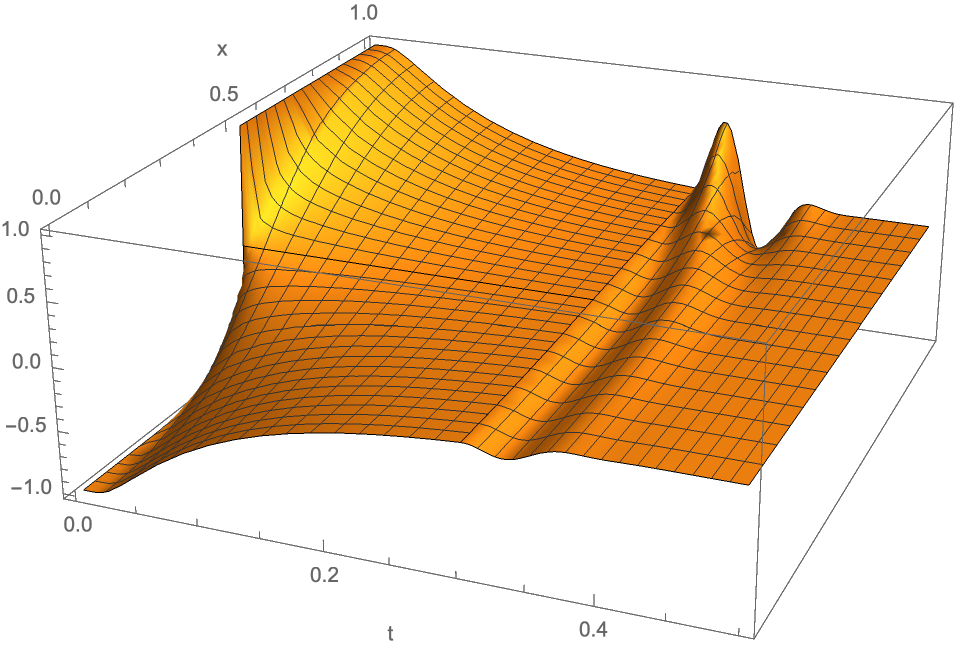}
		\caption{The controlled solution $u(x,t)$  for $x\in(0,1)$ and $t\in(0,1/2)$.}
		\label{fig4}
	\end{figure}
	
	%
	
	%

	\subsection{Discussion on the optimal values for $\tau$ and $N$}\label{subsecopt}
	
	{ We emphasize the fact that this algorithm does not need any fine tuning of the parameters $\tau$ and $N$. This robustness is inherited mainly by the effectiveness of the integral representation of the solution given by the Fokas method. In what follows we discuss the results presented in Tables \ref{tab1} and \ref{tab2} in order to provide an insight to the interested reader who would attempt to realise some experimental simulation.
	}
	
	We perform this algorithmic procedure for different $N=4,...,10$ and $\tau=0, 0.1, 0.15, 0.2, 0.3, 0.35$ and we find $u(x,T)\equiv 0$ up to an \textit{exponentially} small error on $N$ for all $\tau$; the results for $L^2\left[u(x,1/2)\right]$ are presented in Table \ref{tab1}.
	
	\begin{table}[h]
	\resizebox{.7\textwidth}{!}{%
			\begin{tabular}{|l|l|l|l|l|l|l|}
				\hline
				N/$\tau$ & 0        & 0.1      & 0.15     & 0.2      & 0.3      & 0.35     \\ \hline
				4   & 2.54E-08 & 9.91E-08 & 2.32E-07 & 6.40E-07 & 1.00E-05 & 7.26E-05 \\ \hline
				5   & 2.68E-10 & 1.56E-09 & 4.66E-09 & 1.70E-08 & 5.58E-07 & 6.83E-06 \\ \hline
				6   & 2.13E-12 & 1.86E-11 & 7.12E-11 & 3.45E-10 & 2.40E-08 & 4.98E-07 \\ \hline
				7   & 1.33E-14 & 1.76E-13 & 8.61E-13 & 5.57E-12 & 8.25E-10 & 2.91E-08 \\ \hline
				8   & 7.88E-17 & 1.35E-15 & 8.47E-15 & 7.31E-14 & 2.32E-11 & 1.39E-09 \\ \hline
				9   & 3.28E-19 & 1.00E-17 & 8.11E-17 & 7.98E-16 & 5.43E-13 & 5.58E-11 \\ \hline
				10  & 1.03E-21 & 4.80E-20 & 5.01E-19 & 7.75E-18 & 1.08E-14 & 1.90E-12 \\ \hline
			\end{tabular}
		}
		\caption{\label{tab1}$L^2(0,1)$ norm of  $u(x,1/2)$ for different values of $(N,\tau)$.}
	\end{table}
	
	\begin{table}[h]
	\resizebox{.7\textwidth}{!}{%
			\begin{tabular}{|l|l|l|l|l|l|l|}
				\hline
				N/$\tau$ & 0        & 0.1      & 0.15     & 0.2      & 0.3      & 0.35     \\ \hline
				4   & 0.324965 & 0.271611 & 0.260814 & 0.264134 & 0.365895 & 0.59269  \\ \hline
				5   & 0.453611 & 0.370443 & 0.353265 & 0.357046 & 0.506333 & 0.855498 \\ \hline
				6   & 0.596564 & 0.479487 & 0.455493 & 0.460576 & 0.669628 & 1.174559 \\ \hline
				7   & 0.75259  & 0.598281 & 0.567309 & 0.574809 & 0.857294 & 1.557974 \\ \hline
				8   & 0.920823 & 0.726547 & 0.688652 & 0.699917 & 1.070886 & 2.011928 \\ \hline
				9   & 1.100648 & 0.864125 & 0.819538 & 0.836123 & 1.312071 & 2.543485 \\ \hline
				10  & 1.29162  & 1.010936 & 0.960037 & 0.983685 & 1.582591 & 3.160198 \\ \hline
			\end{tabular}
		}
		\caption{\label{tab2}$L^2$ norm of control $h(t)$ for different values of $(N,\tau)$.}
	\end{table}

	The analogue results for the $L^2$ norm of the control, namely $\left(\int_0^T  h^2(s)ds\right)^{1/2}, \ $ are presented in Table \ref{tab2} with respect to $\tau$ and $N$.

	We observe that the values are relatively small, and that we obtain the smallest ones for $\tau=0.15$.  These two Tables also show that for $\tau>0.15$ both the norm of the error and the norm of the control increase. For the value $\tau=0.15$ the error is smaller than $10^{-2(N-1)}$, for all $N=4,...,10$. Furthermore, the norm of the control is increasing with $N$, but is always smaller than 1 for all  $N=4,...,10$; see Table \ref{tab2}.
	
	%
	%
	%
	
	The outcome is that if the control acts immediately, namely $\tau=0$, then the accuracy is better, but then the control is slightly larger. In the next examples we present the results for  $N=6$, and we choose $\tau=0$ and $\tau=0.15$, which yield errors $2 \times 10^{-12}$ and $7 \times10^{-11}$, respectively. These results are presented in figures \ref{fig8}-\ref{fig11}, namely the controls $h(t)$ and the controlled solutions $u(x,t)$.

	%

	
	%
	%

	%

	\begin{figure}[!htbp]
		
		\centering
		\noindent\includegraphics[width=0.6\linewidth ]{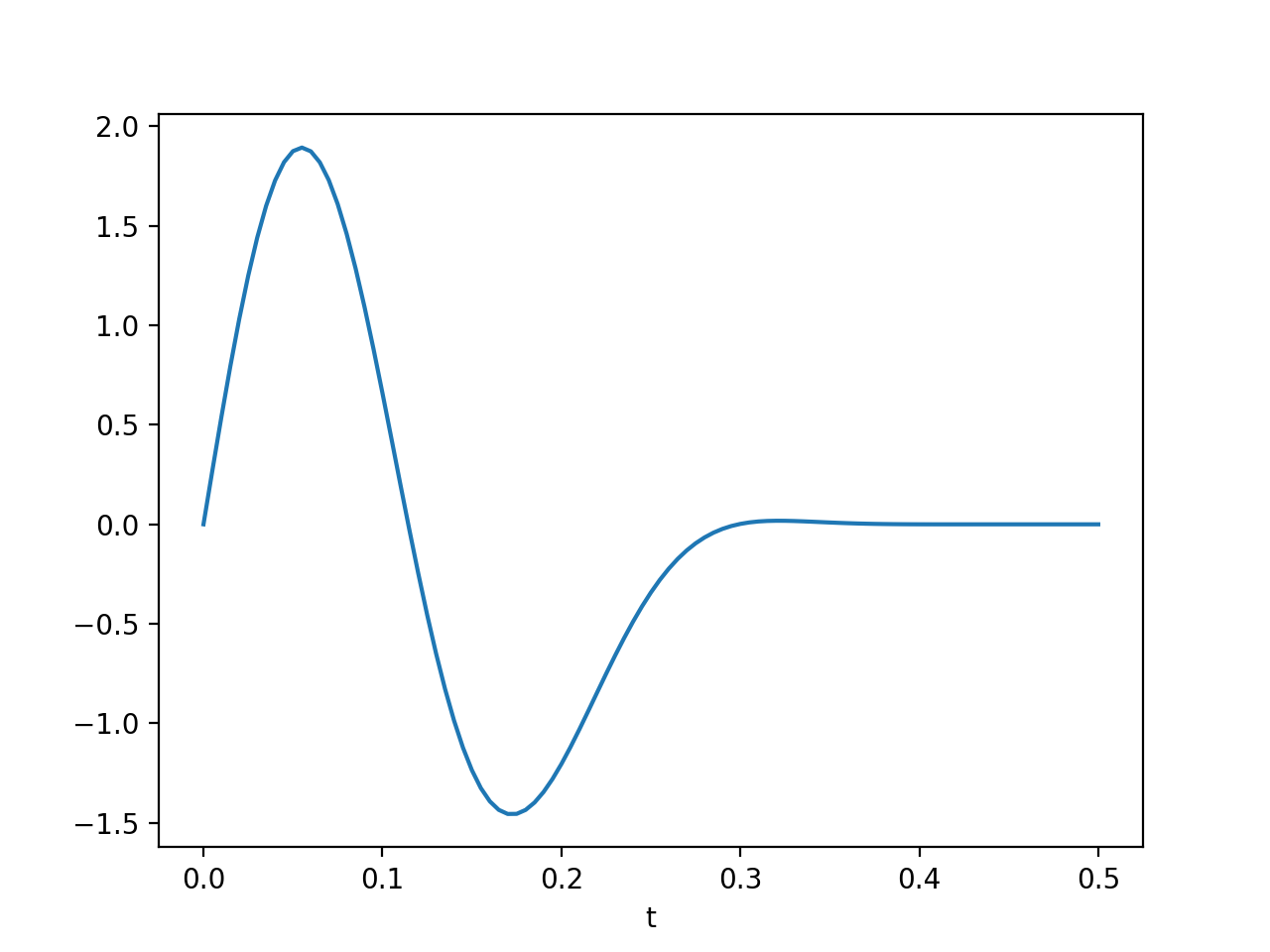}
		\caption{The control  $h(t)$  for $N=6$ and $\tau=0$.}
		\label{fig8}
	\end{figure}
	\begin{figure}
		
		\centering
		\noindent\includegraphics[width=0.6\linewidth ]{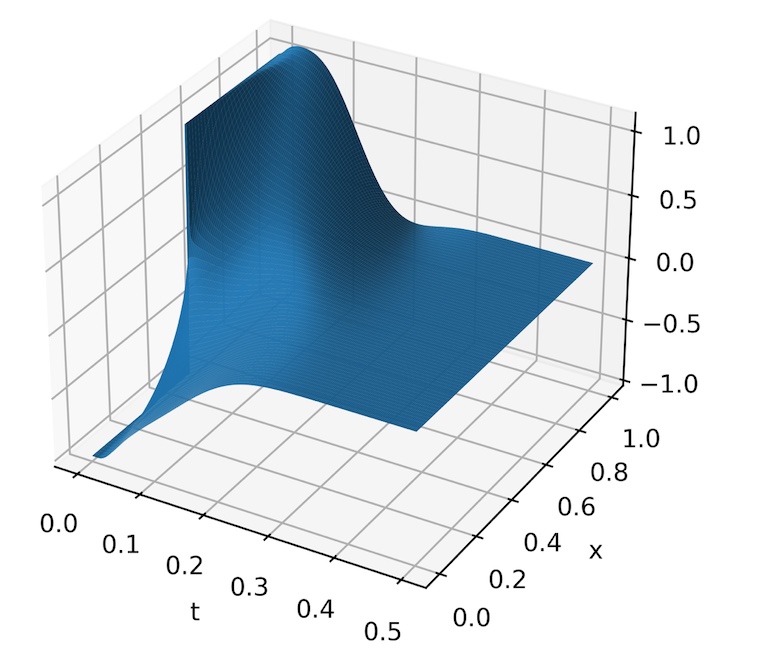}
		\caption{The controlled solution $u(x,t)$  for $x\in(0,1)$ and $t\in(0,1/2)$, for $N=6$ and $\tau=0$.}
		\label{fig9}
	\end{figure}
	\begin{figure}[!htbp]
		
		\centering
		\noindent\includegraphics[width=0.6\linewidth ]{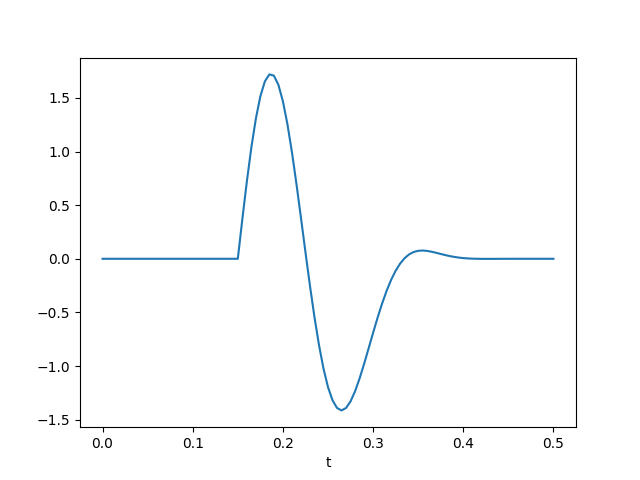}
		\caption{The control  $h(t)$  for $N=6$ and $\tau=0.15$.}
		\label{fig10}
	\end{figure}
	\begin{figure}[!htbp]
		
		\centering
		\noindent\includegraphics[width=0.6\linewidth ]{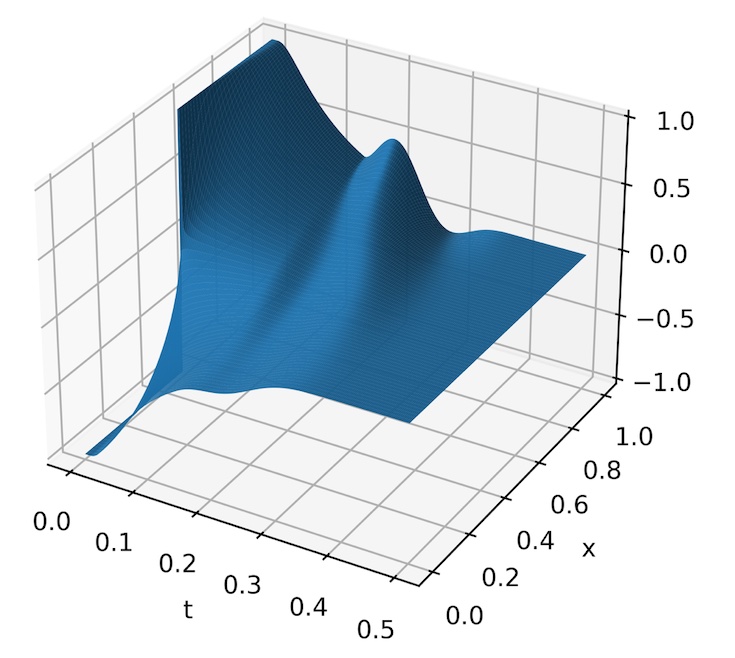}
		\caption{The controlled solution $u(x,t)$  for $x\in(0,1)$ and $t\in(0,1/2)$, for $N=6$ and $\tau=0.15$.}
		\label{fig11}
	\end{figure}
	
	\subsection{Discussion on the discretisation}
	
	{The distribution of $\{x_k\}_{k=0}^{N}$ seems to affect only slightly the performance of the algorithm. However, since the control $h(t)$ acts on the endpoint $x=1$, we expect to have slightly better results if we did not choose naively a uniform distribution, but instead we have chosen a distribution of $\{x_k\}_{k=0}^{N}$ which is denser near $x=1$ and sparser near $x=0$. For matters of completeness we study the example where} the sample is distributed as $x_k=1-\left(\frac{k}{N}\right)^{3/2}$. In this case Table \ref{tab3} demonstrates values of the error up to 6 times smaller than the respective values of Table \ref{tab1}. Furthermore, the control is slightly smaller, comparing the results of Table \ref{tab4} with the ones of Table \ref{tab2}.

	\begin{table}[h]
		\resizebox{.7\textwidth}{!}{%
			\begin{tabular}{|l|l|l|l|l|l|l|}
				\hline
				N/$\tau$ & 0        & 0.1      & 0.15     & 0.2      & 0.3      & 0.35     \\ \hline
				4   & 2.34E-08 & 8.89E-08 & 2.04E-07 & 5.44E-07 & 7.39E-06 & 4.49E-05 \\ \hline
				5   & 2.13E-10 & 1.21E-09 & 3.53E-09 & 1.25E-08 & 3.57E-07 & 3.65E-06 \\ \hline
				6   & 1.38E-12 & 1.18E-11 & 4.43E-11 & 2.08E-10 & 1.27E-08 & 2.20E-07 \\ \hline
				7   & 6.79E-15 & 8.77E-14 & 4.22E-13 & 2.65E-12 & 3.45E-10 & 1.02E-08 \\ \hline
				8   & 2.09E-17 & 5.11E-16 & 3.16E-15 & 2.66E-14 & 7.44E-12 & 3.80E-10 \\ \hline
				9   & 6.40E-20 & 1.92E-18 & 1.53E-17 & 1.73E-16 & 1.31E-13 & 1.15E-11 \\ \hline
				10  & 2.04E-22 & 9.33E-21 & 9.60E-20 & 1.45E-18 & 1.90E-15 & 2.89E-13 \\ \hline
			\end{tabular}
		}
		\caption{\label{tab3}$L^2(0,1)$ norm of  $u(x,1/2)$ for different values of $(N,\tau)$.}
	\end{table}
	
	\begin{table}[h]
		\resizebox{.7\textwidth}{!}{%
			\begin{tabular}{|l|l|l|l|l|l|l|}
				\hline
				N/$\tau$ & 0        & 0.1      & 0.15     & 0.2      & 0.3      & 0.35     \\ \hline
				4   & 0.321827 & 0.26792  & 0.256376 & 0.258205 & 0.347798 & 0.541806 \\ \hline
				5   & 0.450644 & 0.366974 & 0.349079 & 0.351393 & 0.487839 & 0.797046 \\ \hline
				6   & 0.593749 & 0.4762   & 0.451507 & 0.455138 & 0.650918 & 1.111878 \\ \hline
				7   & 0.749906 & 0.595141 & 0.56348  & 0.569531 & 0.838297 & 1.491028 \\ \hline
				8   & 0.918252 & 0.723526 & 0.684946 & 0.694756 & 1.051528 & 1.94067  \\ \hline
				9   & 1.098174 & 0.861202 & 0.81593  & 0.831047 & 1.292284 & 2.467745 \\ \hline
				10  & 1.28923  & 1.008095 & 0.956506 & 0.97867  & 1.562318 & 3.079762 \\ \hline
			\end{tabular}
		}
		\caption{\label{tab4}$L^2$ norm of control $h(t)$ for different values of $(N,\tau)$.}
	\end{table}

	\subsection{Numerical performance}\label{subsecpht}
	
	The numerical calculations were performed in Python using the mpmath \cite{J13} library, which is suitable for real and complex floating-point arithmetic with arbitrary precision. The calculations run with a precision of 30 digits. 	 { Here, we present results for $N\le 10$. The same trend holds also for larger values of $N$, but one has to increase the precision on the calculations.}
	
	Some of the results were obtained also by Mathematica, using the basic commands `NIntegrate', `Inverse' for the numerical integration and the matrix inversion, respectively. The results were identical with the ones obtained by Python, and the computational cost was essentially the same, namely few seconds for the numerical derivation of the control $h(t)$.
	
	\section{Towards a more general set up} \label{sectow}
	{ In this section, we discuss how our algorithm can be modified to treat more general problems in terms of initial-boundary data and control time. A more systematic study of these will be presented in a future work.}

	\subsection{Second specific example} \label{subsecsecond}
	
	The case of initial datum
	$$ u_0(x)=-\cos\left(\frac{\pi}{L} x\right) ,  \qquad x \in (0,L)$$
	yields the RHS of \eqref{int-eq-def}
	$$ P(x;L,T)=-\pi \cos\left(\frac{\pi}{L} x\right) e^{-\left(\frac{\pi}{L}\right)^2 T}, \qquad x \in (0,L).$$
	Indeed, for this initial datum we obtain
	\begin{equation}\label{hat-u-1}
		\hat{u}_0(\lambda)=\frac{i \lambda }{\lambda ^2 -\frac{\pi ^2}{L^2}} +\frac{i \lambda  e^{-i \lambda  L}}{\lambda ^2 -\frac{\pi ^2}{L^2}}, \qquad \lambda\in \C.
	\end{equation}
	By employing in the RHS of \eqref{def:P} we obtain
	\begin{align*}
		P(x;L,T)&=\frac{1}{2} \int_{-\infty}^{\infty} e^{i\lambda x - \lambda^2 T} \left( \frac{i \lambda }{\lambda ^2 -\frac{\pi ^2}{L^2}} +\frac{i \lambda  e^{-i \lambda  L}}{\lambda ^2 -\frac{\pi ^2}{L^2}} \right) d\lambda\\
		&+\frac{1}{2} \int_{\partial D^+}e^{i\lambda x - \lambda^2 T} \frac{i \lambda }{\lambda ^2 -\frac{\pi ^2}{L^2}}   d\lambda-\frac{1}{2} \int_{\partial D^-}e^{i\lambda x - \lambda^2 T}  \frac{i \lambda  e^{-i \lambda  L}}{\lambda ^2 -\frac{\pi ^2}{L^2}}  d\lambda.
	\end{align*}
	{The first integral has two removable singularities at $\lambda=\pm\frac{\pi}{L}$, hence we can deform the real line of integration at these two points by two small  semicircles, which belong at the lower imaginary $\lambda$-plane.  Then, the first term of the first integral of the above equation can be deformed to the $\partial D^+$ yielding two residue contributions at the poles $\lambda=\pm\frac{\pi}{L}$. The second term of the first integral  can be deformed to the $\partial D^-$ yielding zero contribution. Then, Cauchy's theorem yield
		\begin{align*}
			&P(x;L,T)=\frac{1}{2} 2\pi i \text{Res}\left\{ e^{i\lambda x - \lambda^2 T} \frac{i \lambda }{\lambda ^2 -\frac{\pi ^2}{L^2}} \right\}_{\lambda=\pm\frac{\pi}{L}}\\
			&=-\frac{\pi}{2} e^{i\frac{\pi}{L} x - \left(\frac{\pi}{L}\right)^2 T}-\frac{\pi}{2} e^{-i\frac{\pi}{L} x -\left(\frac{\pi}{L}\right)^2 T}\\
			&=-\pi \cos\left(\frac{\pi}{L} x\right) e^{-\left(\frac{\pi}{L}\right)^2 T}.
	\end{align*}}
	

	Fixing $L=1$ and $T=1/2$, then  $$ G(x)= -\pi  \cos\left(\pi x\right) e^{-\frac{\pi^2}{2}}, \ \  x \in (0,1).$$
	For the numerical implementation we need to compute \textit{only} the vector $\mathfrak{g}=\big( G(x_1), \ldots, G(x_{N+1})\big)$, whereas we use the matrix $\mathcal{F}$ obtained already in section \ref{secimp}. Namely, this modification affects only the third step of the numerical implementation of section \ref{subsecnum}.
	
	The figures \ref{fig12}-\ref{fig15} present the analogue of the figures \ref{fig8}-\ref{fig11}, respectively: These are the controls $h(t)$ and the controlled solutions $u(x,t)$ for  $N=6$, and we choose $\tau=0$ and $\tau=0.15$; they yield errors $2 \times 10^{-12}$ and $7 \times10^{-11}$, respectively.
	
	\begin{figure}[!htbp]
		\centering
		\noindent\includegraphics[width=0.6\linewidth ]{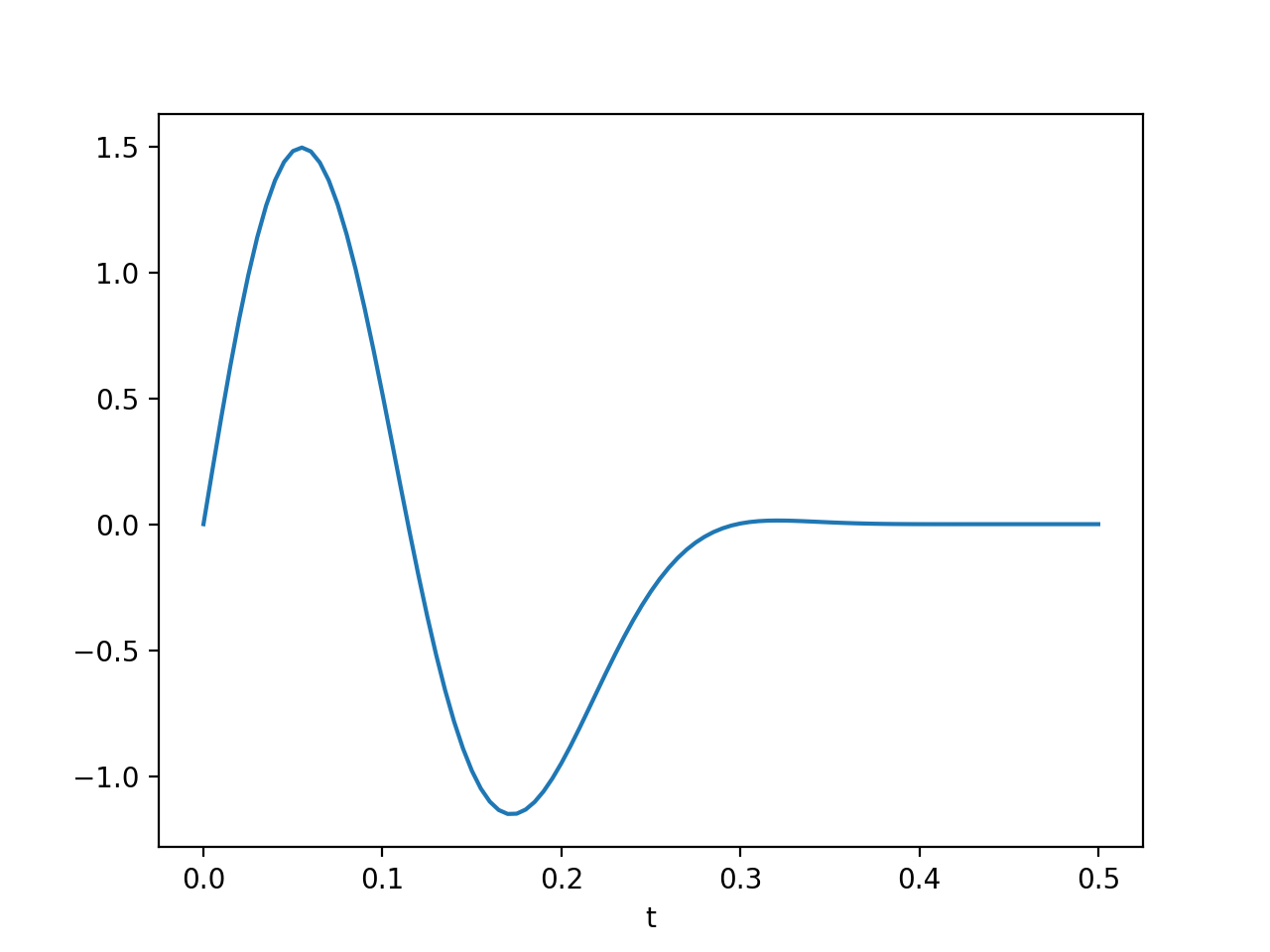}
		\caption{The control  $h(t)$  for $N=6$ and $\tau=0$.}
		\label{fig12}
	\end{figure}
	
	\begin{figure}[!htbp]
		\centering
		\noindent\includegraphics[width=0.6\linewidth ]{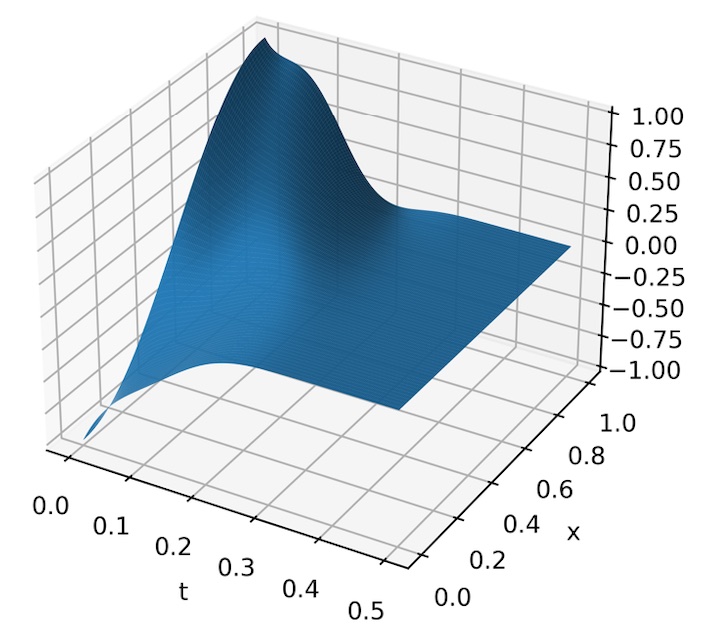}
		\caption{The controlled solution $u(x,t)$  for $x\in(0,1)$ and $t\in(0,1/2)$, for $N=6$ and $\tau=0$.}
		\label{fig13}
	\end{figure}
	
	\begin{figure}[!htbp]
		\centering
		\noindent\includegraphics[width=0.6\linewidth ]{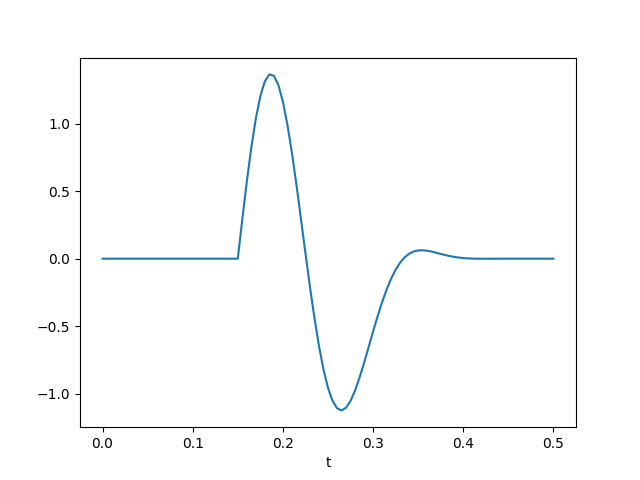}
		\caption{The control  $h(t)$  for $N=6$ and $\tau=0.15$.}
		\label{fig14}
	\end{figure}
	
	\begin{figure}[!htbp]
		\centering
		\noindent\includegraphics[width=0.6\linewidth ]{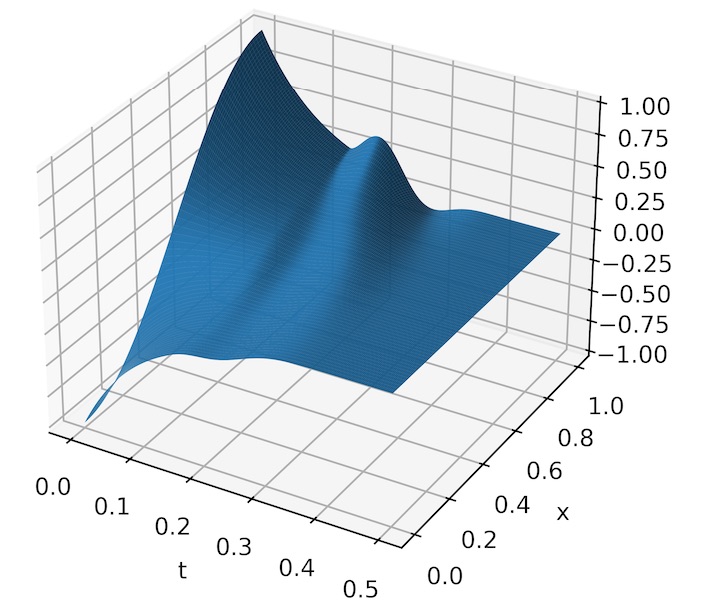}
		\caption{The controlled solution $u(x,t)$  for $x\in(0,1)$ and $t\in(0,1/2)$, for $N=6$ and $\tau=0.15$.}
		\label{fig15}
	\end{figure}

	\subsection{The general basis}\label{subsecgen}
	
	It is well known that a large class of  continuous initial data can be approximated by the cosine-Fourier series, thus we find useful to discuss the case of initial value
	$$ u_0^{(m)}(x)=-\cos\left(m\frac{\pi}{L} x\right) ,  x \in (0,L), \ \ m=1,2,\ldots.$$
	This yields the RHS of \eqref{int-eq-def}
	$$ P^{(m)}(x;L,T)=-\pi\cos\left(m\frac{\pi}{L} x\right) e^{-\left(m\frac{\pi}{L}\right)^2 T}, \quad x \in (0,L),\ \ m=1,2,\ldots,$$
	which gives for $L=1$ and $T=1/2$
	$$ G^{(m)}(x)=-\pi\cos\left(m \pi x\right) e^{-\frac{m^2 \pi^2}{2}}, \quad x \in (0,1),\ \ m=1,2,\ldots.$$
	The only modification of the algorithm presented in section \ref{subsecnum} occurs on step [3.], namely updating the RHS of \eqref{int-eq-g-num-2}, whereas the matrix $\mathcal{F}$ remains the same as in section \ref{secimp}.

	\subsection{General initial conditions}
	
	The above result indicates a pathway for treating the problem fully numerically for a general initial condition $u_0(x)$, under the requirement that one can effectively approximate  $u_0(x)$ by its cosine-Fourier series in the interval $(0,L)$. 
	{ The concept is rather straightforward: First, write $u_0(x)$ as (truncated) cosine-Fourier series: $$u_0(x)\approx\frac{c_0}{2}+\sum_{m=1}^M c_m u_0^{(m)}(x),$$ with $u_0^{(m)}(x)$ given in subsection \ref{subsecgen}; namely, we obtain the real constants $$c_m = -\frac{2}{L} \int_0^L u_0(x) \cos \left(\frac{m\pi }{L} x\right)  dx.$$\\}\\
	{ Second, employ the procedure of section \ref{subsecsecond} with 	$$P(x;L,T) \approx\frac{c_0}{2} \pi+\sum_{m=0}^M c_m P^{(m)}(x;L,T),$$ where $P^{(m)}(x;L,T)$ are given in subsection \ref{subsecgen}.  We note that this procedure may be quicker than the one in section \ref{secfirst} because (since we have obtained $c_m$) the function $G(x)=P(x;L,T)$ in the RHS of \eqref{int-eq-g-num-2} is available analytically, in contrast to the numerical integration needed in section \ref{secfirst}; furthermore, the fact that $P^{(m)}(x;L,T)$ is exponentially decaying with $m$, suggests that one only needs very few values of $c_m$ to approximate $P(x;L,T)$, up to exponentially small error.\\}
	{ In fact, if one applies this approach to the example of section \ref{secfirst}, using only the first term of the cosine series of the initial data $u_0(x)$, then the function $G(x)$ of \eqref{G-example1}  will be approximated by $\tilde{G}(x)= -4 \cos(\pi x) e^{-\pi^2/2}$. A simple numerical calculation yields that $|G-\tilde{G}|$ is $O\left(10^{-20}\right)$.}
	
	Also, one could appropriately modify the analysis of the section \ref{secfirst} in order to incorporate the Haar-wavelet basis. Both bases as well as different type of boundary controls, e.g., Dirichlet and Robin boundary conditions will be analysed in a future work.
	
	
	\subsection{External boundary sources and control}
	Most well-known boundary control methods in the literature assume homogeneous boundary conditions on the part of boundary where control is not supplied.  Although such an assumption simplifies the mathematical challenges, it is not always realistic. The properties of an evolution on the part of boundary layer may vary with respect to time and (in higher dimensions {also} with respect to) space due to external {sources}.  {A physical example is the evolution of temperature on a finite rod which does not satisfy ideal/perfect insulation properties at one of its ends where an external heat source can manipulate the evolution.  A trivial case is a constant external heat being supplied at one end of the rod while a control is acting at the other end. This could for instance occur if the rod is in contact with another - say relatively larger material/medium that supplies a constant flow of heat from one of ends of the rod. Of course, here we are assuming that this external source can be characterized as a deterministic data. Under this assumption, it is desirable to introduce a control method which works in presence of external boundary sources.} Considering the system \eqref{heat-fi}, the input $g(t)$ can be used to denote an external {source} while $h(t)$ denotes a sought-after control (or vice versa). The effect of such an external input} in our algorithm would only give rise to one additional integral that involves $g$ at the right hand side of \eqref{def:P}.  Our algorithm easily adapts to this modification. 
{Clearly, an open loop control system as in this paper will not make much sense for stochastic disturbances and noises which are not known a priori.  Such noises should be treated instead with closed loop control systems/stabilizers.}  To the best of our knowledge, other control methods do not easily adapt to external {source} case. {For instance, it is not clear to us how a control method will adapt to external source case if it uses smoothing and intrinsic decay properties of heat equation as these features may no longer be valid in presence of boundary data}.  This rather strong regularity property is no longer true if there is an external rough input at one end of the boundary even if no control is applied at the other end.
{It should be noted however that in the case of an external source even if the sought-after control steers solutions to zero at time $t=T$, the solution will in general not stay at the zero state for $t>T$.  However, if the external manipulation has some rather nice dynamical character such as being an exponentially decaying function, then one can anticipate that the solution will remain near zero.  Of course, in the case the external manipulation is abrupt, in addition to having a large control effort, there may also be the situation that solution will quickly diverge away from zero state for $t>T$. Clearly, the case of an external manipulation may lead to numerous physical situations. A more detailed analysis of such interesting phenomena will be presented in a future work.}

\subsection{Control in small time}
Note that Theorem \ref{controlthm} does not put a restriction on the size of $T$. Namely, null controllability holds no matter how small or large $T$ is. This is aligned with the infinite speed of propagation of solutions of the heat equation. Null controllability of heat for arbitrarily small $T$ contrasts with the case of hyperbolic PDEs such as the wave equation. For the latter, there is a lower bound for $T$ in order for controllability to hold, which reflects the fact that solutions propagate with finite speed. Although, heat is null controllable in arbitrarily small time, the control effort ($L^2(0,T)$-norm of $h$) will get larger as $T$ gets smaller.  Preliminary results of our algorithm suggest a sub-exponential increase (with respect to $T$) of the control $h(t)$ as $T$ approaches $0$.

{	
	Indeed, for the case of the example of section \ref{secfirst}, with $\tau=0$, our algorithm yields the following results:
	\begin{itemize} 
		\item For $T=0.1$, and $N=10$, we obtain a control $h$ whose magnitude is of order $10^3$.
		\item For $T=0.01$, and $N=48$, we obtain a control $h$ whose magnitude is of order $10^{24}$. 
	\end{itemize}
	In both cases $u(x,T)$ vanishes up to an error of order $10^{-8}$.  The rigorous asymptotic analysis as $T \to 0$, both numerically and analytically will be performed in future work. }

\section{Conclusion and discussion of the results}
With this work we aim to complete the introduction of a new methodology, based on the Fokas method, for treating problems originated in the control theory. In \cite{KO20} we proved the lack of null controllability of the heat equation on half line (and in arbitrarily many spatial dimensions). In the current work we presented in detail a numerical algorithm to obtain a null control of Neumann type for the heat equation on the finite interval. The efficiency of the algorithm, which reflects the numerical power of the Fokas method, is illustrated in section \ref{secimp}, considering both its low computational cost and its high accuracy.
The rather straightforward numerical algorithm aims to solve equation \eqref{int-eq-def} which follows from equation \eqref{sol-1}, namely the integral representation of the solution which is provided by the Fokas method. The unified character of the method reduces the complexity of the derivation of the integral representation to an exercise in undergraduate books for a wide range of initial boundary value problems for linear PDEs, see \cite{FKinpress}; already the analogue of \eqref{sol-1} is readily available for many different problems in the vast literature. In analogy, the other basic element of the Fokas method, namely the global relation, is the basis for analyzing theoretical aspects of controllability, \cite{KO20}. The afore-mentioned results provide the background for a coherent and easy-to-apply methodology for treating control problems both theoretically and computationally.

\end{document}